\documentstyle[aps,pra,amsfonts]{revtex}
\begin{document} 

\begin{titlepage}
\title{A Note on Regularized Shannon's Sampling Formulae} 
\author{ Liwen Qian and G. W. Wei\\ 
\emph {Department of Computational Science, 
National University of Singapore\\
Singapore 117543}}

\date{\today} 
\maketitle

\begin{abstract} 
Error estimation is given for a 
regularized Shannon's sampling formulae, 
which was found 
to be accurate and robust for numerically 
solving partial differential equations.
\end{abstract} 



\end{titlepage}

\section{Introduction} 

In previous work\cite{weijcp}, one of the present authors
proposed a discrete singular convolution (DSC) algorithm  for 
computer realization of singular convolutions involving 
singular kernels of delta type, Abel type and Hilbert type.
One of illustrations for the algorithm was 
Shannon's sampling formulae\cite{Shannon} which 
plays an important 
role in the approximation of the delta distribution and 
generalized derivatives\cite{Chui}.
However, in practical computations, the truncation error 
of Shannon's sampling formulae is substantial\cite{Marks}.
A regularization technique\cite{weiprl} was used to  
construct a regularized Shannon's sampling 
formulae\cite{weijcp,weicpl}, which was found to be 
extremely accurate and robust for resolving various 
challenging dynamical problems, such as the homoclinic 
orbit excitation of the Sine-Gordon equation\cite{weiphysica}, 
the Navier-Stokes flow in complex geometries\cite{wanwei}, 
shock capturing of the inviscid Burgers' equation\cite{weigu},
molecular quantum  system described by the 
Schr\"{o}dinger equation\cite{weijpb}
and nonlinear pattern formation of 
the Cahn-Hilliard equation\cite{GLW}. 

The objective of the present note is to provide
a theoretical analysis for the previous excellent 
numerical results. Rigorous    
error estimations of the regularized 
Shannon's sampling formulae are given for 
their applications  
to interpolations and derivatives of a function.

\section{Main result} 
{\bf Theorem.} Let $f$ be a function $f\in L^2(R)\cap C^{s}(R)$ 
and bandlimited to 
$B$, $(B< \frac{\pi}{\Delta},~ \Delta$ is the grid spacing).  
For a fixed $t\in R $ 
and $\sigma>0$, denote $g(x)=f(x)H_k(\frac{t-x}{\sqrt{2}\sigma})$, 
where $H_k(x)$ is the
$k$th order Hermite polynomial. If $g(x)$ satisfies 
\begin{equation} 
g'(x)\leq
g(x)\frac{(x-t)}{{\sigma}^2} 
\end{equation} 
for $x\geq t+(M_1-1)\Delta$, and
\begin{equation} 
g'(x)\geq g(x)\frac{(x-t)}{{\sigma}^2} 
\end{equation} 
for $x\leq
t-M_2\Delta$, where $M_1,M_2\in \mathcal{N}$,  then for any 
$s\in \mathcal{Z}^{+}$
\begin{eqnarray} 
&&\left\|
f^{(s)}(t)-\sum_{n=\lceil\frac{t}{\Delta}\rceil-M_2}^{\lceil
\frac{t}{\Delta}
\rceil+M_1}f(n\Delta)
\left[\frac{\sin\frac{\pi}{\Delta}(t-n\Delta)}{\frac{\pi}{\Delta}
(t-n\Delta)}
\exp(-\frac{(t-n\Delta)^2}{2{\sigma}^2})\right]^{(s)}
\right\|_{L^2(R)}\nonumber\\ 
&&\leq
\sqrt{3}\left[\frac{\| f^{(s)}(t)\|_{L^2(R)}}{2\pi\sigma
(\frac{\pi}{\Delta}-B)
\exp(\frac{{\sigma}^2(\frac{\pi}{\Delta}-B)^2}{2})}\right.
\nonumber\\ 
&&\left.+ \frac{\|
f(t)\|_{L^2(R)}\sum_{i+j+k=s}\frac{s!{\pi}^{i-1}
H_k(\frac{-M_1\Delta}{\sqrt{2}\sigma})}
{i!k!{\Delta}^{i-1}(\sqrt{2}\sigma)^k((M_1-1)\Delta)^{j+1}}}
{\exp(\frac{(M_1\Delta)^2}{2{\sigma}^2})} \right.\nonumber\\ 
&&\left.+\frac{\|
f(t)\|_{L^2(R)}\sum_{i+j+k=s}\frac{s!{\pi}^{i-1}
H_k(\frac{-M_2\Delta}{\sqrt{2}\sigma})}
{i!k!{\Delta}^{i-1}(\sqrt{2}\sigma)^k(M_2\Delta)^{j+1}}}
{\exp(\frac{(M_2\Delta)^2}{2{\sigma}^2})}\right],
\label{eq3} 
\end{eqnarray}
where superscript, $(s)$, denotes the $s$th order derivative.

\section{Proof}

\subsection{Separation of the error}

The error breaks  naturally into a few components. 
Denote 
\begin{eqnarray}
&&E(t)=f^{(s)}(t)-\sum_{n=\lceil\frac{t}{\Delta}\rceil-M_2}^{n=\lceil
\frac{t}{\Delta}
\rceil+M_1}f(n\Delta)
\left[\frac{\sin\left(\frac{\pi}{\Delta}(t-n\Delta)\right)}{\frac{\pi}
{\Delta}(t-n\Delta)}
\exp\left(-\frac{(t-n\Delta)^2}{2{\sigma}^2}\right)\right]^{(s)}\\
&&E_1(t)=\sum_{n=-\infty}^{n=+\infty}f(n\Delta)\left[\frac{\sin
\left(\frac{\pi}
{\Delta}(t-n\Delta)\right)}
{\frac{\pi}{\Delta}(t-n\Delta)}
-\frac{\sin\left(\frac{\pi}{\Delta}(t-n\Delta)\right)}{\frac{\pi}
{\Delta}(t-n\Delta)}
\exp\left(-\frac{(t-n\Delta)^2}{2{\sigma}^2}\right)\right]^{(s)}
\label{eq5}
\\
&&E_2(t)=\sum_{n\geq\lceil\frac{t}{\Delta}\rceil+M_1}f(n\Delta)\left[
\frac{\sin\left(\frac{\pi}
{\Delta}(t-n\Delta)\right)} {\frac{\pi}{\Delta}(t-n\Delta)}
\exp\left(-\frac{(t-n\Delta)^2}{2{\sigma}^2}\right)\right]^{(s)}
\label{eq6}
\\
&&E_3(t)=\sum_{n\leq\lceil\frac{t}{\Delta}\rceil-M_2}f(n\Delta)\left[
\frac{\sin\left(
\frac{\pi}
{\Delta}(t-n\Delta)\right)} {\frac{\pi}{\Delta}(t-n\Delta)}
\exp\left(-\frac{(t-n\Delta)^2}{2{\sigma}^2}\right)\right]^{(s)}. 
\label{eq7}
\end{eqnarray} 
Here, $E_1(t)$ is regularization error. 
$E_2(t)$  and $E_3(t)$  are truncation errors.
From Shannon's sampling theorem \cite{Shannon} 
\begin{equation}
f(t)=\sum_{n=-\infty}^{n=+\infty}f(n\Delta)\frac{\sin
\left(\frac{\pi}{\Delta}
(t-n\Delta)\right)}
{\frac{\pi}{\Delta}(t-n\Delta)}, 
\end{equation} 
which can be differentiated term by term,  
the total error can be written as a sum of three components 
\begin{equation} 
E(t)=\sum_{i=1}^{3}E_i(t).
\end{equation} 
The corresponding error norms satisfy
\begin{equation} \label{eq10} 
\| E(t)\|_{L^2(R)}\leq \sqrt{3}\sum_{i=1}^{3}\| E_i(t)\|_{L^2(R)}.
\end{equation}

\subsection{Estimation of $E_1(t)$}

Let $f\hat(\omega)$ be the
	Fourier transform of $f(x)$, and
   $f\hat(\omega)=\int_{R}f(x)\exp(ix\omega)dx$.
	Since 
	\begin{equation}
	\left[\frac{\sin\left(\frac{\pi}{\Delta}(t-n\Delta)\right)}
	{\frac{\pi}{\Delta}(t-n\Delta)}\right]\hat(\omega)=
	\Delta \exp(-in\Delta
	\omega)\chi_{[-\frac{\pi}{\Delta},\frac{\pi}{\Delta}]}(w) 
	\end{equation} 
	and
	\begin{equation}
	\left[\exp\left(-\frac{(t-n\Delta)^2}{2{\sigma}^2}\right)\right]
	\hat(\omega)=\sqrt{2\pi}\sigma
	\exp(-in\Delta \omega- \frac{{\sigma}^2 {\omega}^2}{2}), 
	\end{equation} 
one writes
	\begin{eqnarray}\label{eq13}
	&&\left[\frac{\sin\left(\frac{\pi}{\Delta}(t-n\Delta)\right)}
        {\frac{\pi}
	{\Delta}(t-n\Delta)}\right]\hat(\omega)\ast
	\left[\exp\left(-\frac{(t-n\Delta)^2}{2{\sigma}^2}\right)\right]
        \hat(\omega)
	\nonumber\\
	&&=\int_R \Delta\sqrt{2\pi}\sigma \exp(-in\Delta
	(\omega-\theta))\chi_{[\theta-\frac{\pi}{\Delta},
	\theta+\frac{\pi}{\Delta}]}(\omega) \exp(-in\Delta
	\theta-\frac{{\sigma}^2{\theta}^2}{2})d\theta\nonumber\\
	&&=\Delta\sqrt{2\pi}\sigma \exp(-in\Delta
	\omega)\int_{\theta-\frac{\pi}{\Delta}}^{\theta+\frac{\pi}{\Delta}}
	\exp(-\frac{{\sigma}^2{\theta}^2}{2})d\theta . 
	\end{eqnarray} 
From Eq. (\ref{eq13})  
	\begin{eqnarray}
	&&\left[\sum_{n=-\infty}^{n=+\infty}f(n\Delta)\frac{\sin\left(
        \frac{\pi}{\Delta}
	(t-n\Delta)\right)}{\frac{\pi}
	{\Delta}
	(t-n\Delta)}\exp\left(-\frac{(t-n\Delta)^2}{2{\sigma}^2}\right)\right]
	\hat(\omega)\nonumber\\
	&&=\sum_{n=-\infty}^{n=+\infty}f(n\Delta)
	\frac{1}{2\pi}\left[\frac{\sin\left(\frac{\pi}{\Delta}
	(t-n\Delta)\right)}{\frac{\pi}{\Delta}(t-n\Delta)}]\hat(\omega)\ast
	[\exp\left(-\frac{(t-n\Delta)^2} {2{\sigma}^2}\right)\right]
	\hat(\omega)\nonumber\\
	&&=\sum_{n=-\infty}^{n=+\infty}f(n\Delta)\Delta \exp(-in\Delta
	\omega)\frac{1}{\sqrt{\pi}} \int_{\frac{\sigma
	(\omega-\frac{\pi}{\Delta})}{\sqrt{2}}}^{\frac{\sigma
	(\omega+\frac{\pi}{\Delta})} {\sqrt{2}}}\exp(-t^2)dt. \label{eq14}
	\end{eqnarray} 
Since function $f$ satisfies
	\begin{equation} 
	f\hat(\omega)\in L^2[-B,B]\subset
	L^2[-\frac{\pi}{\Delta},\frac{\pi}{\Delta}], 
	\end{equation} 
it has a Fourier series expansion  
	\begin{equation}\label{eq16} 
	f\hat(\omega)=\sum_{n=-\infty}^{\infty}c_n\exp(in\Delta\omega), 
	\end{equation}
where coefficients is given by 
	\begin{equation}
	c_n=\frac{\Delta}{2\pi}\int_{-\frac{\pi}{\Delta}}^{\frac{\pi}{\Delta}}
	f\hat(\omega)
	\exp(-in\Delta\omega)d\omega=\Delta f(-n\Delta). 
        \end{equation} 
Equivalently, $f\hat(\omega)$ can be written 
	\begin{equation} \label{eq18} 
	f\hat(\omega)=f\hat(\omega)\chi_{[-B,B]}(\omega)
	=\sum_{n=-\infty}^{\infty}\Delta
	f(n\Delta)\exp(-in\Delta\omega)\chi_{[-B,B]}(\omega). 
	\end{equation} 
Denote
	\begin{equation}
	\varepsilon(\omega)=\chi_{[-B,B]}(\omega)-\frac{1}{\sqrt{\pi}}
	\int_{\frac{\sigma (\omega-\frac{\pi}{\Delta})}{\sqrt{2}}}^{\frac{
        \sigma
	(\omega+\frac{\pi}{\Delta})} {\sqrt{2}}}\exp(-t^2)dt, 
	\end{equation} 
then combining Eqs. (\ref{eq5}), (\ref{eq14}), (\ref{eq16}) 
        and (\ref{eq18}),  one has
	\begin{equation}\label{eq20}
	E_1\hat(\omega)=(i\omega)^sf\hat(\omega)\varepsilon(\omega). 
	\end{equation} 
For $\omega\in [-B,B]$, $\varepsilon(\omega)$ can be evaluated as
	\begin{eqnarray}
	\varepsilon(\omega)&=&\frac{1}{\sqrt{\pi}}\left[\int_{-\infty}^{\infty}
	\exp(-t^2)dt-\int_{\frac{\sigma
	(\omega- \frac{\pi}{\Delta})}{\sqrt{2}}}^{\frac{\sigma
	(\omega+\frac{\pi}{\Delta})} {\sqrt{2}}}\exp(-t^2)dt\right]\nonumber\\
	&=&\frac{1}{\sqrt{\pi}}\left[\int_{\frac{\sigma
	(\frac{\pi}{\Delta}-\omega)}{\sqrt{2}}}^{\infty}
	\exp(-t^2)dt+\int_{\frac{\sigma (\omega+\frac{\pi}{\Delta})}
	{\sqrt{2}}}^{\infty}\exp(-t^2)dt\right]. 
	\end{eqnarray} 
Moreover, for $x\geq 0$, the 
        following inequality \cite{Zwillinger} is valid 
	\begin{equation} 
	\frac{1}{x+\sqrt{x^2+2}}\leq
	\exp(x^2)\int_x^{\infty}\exp(-t^2)dt
	\leq \frac{1}{x+\sqrt{x^2+\frac{4}{\pi}}}.
	\end{equation} 
Therefore, the estimation for  $\varepsilon(\omega)$ is obtained as
	\begin{eqnarray} 
	\varepsilon(\omega)&\leq&
	\frac{1}{\sqrt{\pi}}\left(\frac{\exp\left(-\frac{{\sigma}^2
	(\frac{\pi}{\Delta}-\omega)^2}{2}\right)}
	{\sqrt{2}\sigma(\frac{\pi}{\Delta}-\omega)}+\frac{\exp\left(-\frac{{
        \sigma}^2
	(\frac{\pi}{\Delta}+
	\omega)^2}{2}\right)} {\sqrt{2}\sigma(\frac{\pi}{\Delta}+\omega)}
        \right)
	\nonumber\\ &\leq&
	\frac{1}{\sigma(\frac{\pi}{\Delta}-B)\exp\left(\frac{{\sigma}^2(
        \frac{\pi}
	{\Delta}-B)^2}{2}\right)}.    \label{eq23}
	\end{eqnarray} 
It follows from Eqs. (\ref{eq20}) and (\ref{eq23}) that
	\begin{eqnarray} 
	E_1\hat(\omega) \leq
	\frac{f\hat(\omega)(i\omega)^s}{\sigma(\frac{\pi}{\Delta}-B)
	\exp(\frac{{\sigma}^2
	(\frac{\pi}{\Delta}-B)^2}{2})}. 
	\end{eqnarray} 
From the Parseval identity, one has  
	\begin{eqnarray} \label{eq25} 
	\| E_1(t)\|_{L^2(R)}&=&\frac{1}{2\pi}\|
	E_1\hat(\omega)\|_{L^2(R)}\nonumber\\ 
	&\leq&\frac{\|
	f^{(s)}(t)\|_{L^2(R)}}{2\pi\sigma(\frac{\pi}{\Delta}-B)
	\exp\left(\frac{{\sigma}^2(\frac{\pi}{\Delta}-B)^2}{2}\right)}. 
	\end{eqnarray}

\subsection{Estimation of $E_2(t)$}

Differentiations can be written 
\begin{eqnarray}
E_2(t)=\sum_{n\geq \lceil\frac{t}{\Delta}\rceil+M_1}f(n\Delta)
\left[\sum_{i+j+k=s}\frac{s!}{i!k!}(\frac{\pi}{\Delta})^{i-1}
\sin(\frac{\pi}{\Delta}t-n\pi+\frac{\pi i}{2})\right.\nonumber\\
\left.\frac{(-1)^{j}}{(t-n\Delta)^{j+1}}\frac{(-1)^k}{(\sqrt{2}\sigma)^k}
H_k(\frac{t-n\Delta}{\sqrt{2}\sigma})\exp\left(-\frac{(t-n\Delta)^2}
{2{\sigma}^2}\right)\right],
\end{eqnarray} 
where $H_k(x)$ is the Hermite polynomial 
\begin{equation}
\exp(-x^2)H_k(x)=(-1)^k(\frac{d}{dx})^k\exp(-x^2). 
\end{equation} 
Let
$l=n-\lceil\frac{t}{\Delta}\rceil$, where $\lceil x\rceil$ is the 
integral  part of
$x$ and $\lfloor x \rfloor=x-\lceil x\rceil$, then 
\begin{eqnarray}
E_2(t)&=&\sum_{l\geq M_1}f(t+l\Delta-\lfloor\frac{t}{\Delta}\rfloor\Delta)
\left[\sum_{i+j+k=s}\frac{s!}{i!k!}(\frac{\pi}{\Delta})^{i-1}
(-1)^{j+k+l}\right.\nonumber\\
&&\left.\sin\left(\lfloor\frac{t}{\Delta}\rfloor\pi +\frac{\pi i}{2}\right)
\frac{1}{(\sqrt{2}\sigma)^k{\Delta}^{j+1}(-l+\lfloor\frac{t}{\Delta}
\rfloor)^{j+1}}\right.\nonumber\\
&&\left[H_k\left(\frac{-l\Delta+\lfloor\frac{t}{\Delta}\rfloor\Delta}
{\sqrt{2}\sigma}\right)\right]
\exp\left(-\frac{(l\Delta-\lfloor\frac{t}{\Delta}\rfloor\Delta)^2}
{2{\sigma}^2}\right)\nonumber\\
&=&\sum_{i+j+k=s}\sum_{l\geq M_1}F_k(l)s_{i,j,k}(l), 
\end{eqnarray} 
where
\begin{eqnarray}
F_k(l)=&&f(t+l\Delta-\lfloor\frac{t}{\Delta}\rfloor\Delta)H_k\left(
\frac{-l\Delta+
\lfloor\frac{t}{\Delta}\rfloor\Delta}{\sqrt{2}\sigma}\right)\nonumber\\
&&\exp\left(-\frac{(l\Delta-\lfloor\frac{t}{\Delta}\rfloor\Delta)^2}
{2{\sigma}^2}\right)\\
s_{i,j,k}(l)=&&(-1)^l\frac{s!{\pi}^{i-1}(-1)^{j+k}\sin(\lfloor\frac{t}
{\Delta}\rfloor\pi+\frac{\pi
i}{2})}
{i!k!{\Delta}^{i+j}(-l+\lfloor\frac{t}{\Delta}\rfloor)^{j+1}
(\sqrt{2}\sigma)^k}.
\end{eqnarray}

Two simple lemmas are required.

{\bf Lemma 1 (Abel's inequality)}\cite{Mitrinovic}.
 For two sequences $\{a_n\}, 
\{b_n\}, b_1\geq
b_2\geq \dots \geq b_n, a_n, b_n\in R$, set 
\begin{eqnarray} 
s_k&=&\sum_{i=1}^ka_i\\
m&=&\min_{1\leq k\leq n}s_k\\ 
M&=&\max_{1\leq k\leq n}s_k, 
\end{eqnarray} 
then
\begin{equation} 
mb_1\leq \sum_{i=1}^na_ib_i\leq Mb_1. 
\end{equation} 

{\bf Lemma 2}.  As all the notations unchanged, 
set $g(x)=f(x)H_k(\frac{t-x}{\sqrt{2}\sigma})$. 
If $g(x)$
satisfies 
\begin{equation} 
g'(x)\leq g(x)\frac{(x-t)}{{\sigma}^2}, 
\end{equation}
whenever $x\geq t+(M_1-1)\Delta$, then $\{F_k(l)\}_{l\in \mathcal{N}}$ 
is a decreasing sequence.

The proof is obvious by taking the first order derivative.

Let denote
\begin{equation} 
S_{i,j,k}(N)=\sum_{l\geq M_1}^{M_1+N}s_{i,j,k}(l). 
\end{equation}
It is estimated that  
\begin{eqnarray} 
|S_{i,j,k}(N)|&=&\left|\sum_{l\geq
M_1}^{M_1+N}(-1)^l\frac{s!{\pi}^{i-1}\sin(\lfloor\frac{t}{\Delta}
\rfloor\pi +\frac{\pi
i}{2})(-1)^{k+j}} {i!k!{\Delta}^{i+j} (\sqrt{2}\sigma)^k}
\frac{1}{(-l+\lfloor\frac{t}{\Delta}\rfloor)^{j+1}}\right|\nonumber\\
&=&\left|\frac{s!{\pi}^{i-1}\sin(\lfloor\frac{t}{\Delta}\rfloor\pi 
+\frac{\pi
i}{2})(-1)^{k+j}}{i!k!{\Delta}^{i+j} (\sqrt{2}\sigma)^k}
\sum_{l\geq M_1}^{M_1+N}
(-1)^l
\frac{1}{(-l+\lfloor\frac{t}{\Delta}\rfloor)^{j+1}}\right|
\nonumber\\ 
&\leq&
\frac{s!{\pi}^{i-1}}{i!k!{\Delta}^{i+j} (\sqrt{2}\sigma)^k}
\frac{1}{(M_1-1)^{j+1}}.
\end{eqnarray} 
Then from $(28)$ and $(37)$, and by using lemma $1$ and lemma $2$,
one has 
\begin{eqnarray}
E_2(t)\leq
f(t+M_1\Delta)\exp\left(-\frac{(M_1\Delta)^2}{2{\sigma}^2}\right)
\nonumber\\
\times\sum_{i+j+k=s}\frac{s!{\pi}^{i-1}H_k(\frac{-M_1\Delta}{\sqrt{2}
\sigma})}
{i!k!{\Delta}^{i-1}(\sqrt{2}\sigma)^k\left((M_1-1)\Delta\right)^{j+1}}. 
\end{eqnarray}
This gives rise to 
\begin{equation}\label{eq39}  
\| E_2(t)\|_{L^2(R)}\leq \frac{\| f(t)\|_{L^2(R)}
\sum_{i+j+k=s}\frac{s!{\pi}^{i-1}H_k(\frac{-M_1\Delta}
{\sqrt{2}\sigma})}{i!k!{\Delta}^{i-1}(\sqrt{2}\sigma)^k
\left((M_1-1)\Delta\right)^{j+1}}}
{\exp\left(\frac{(M_1\Delta)^2}{2{\sigma}^2}\right)}. 
\end{equation}

\subsection{Estimation of  $E_3(t)$} 

A result like lemma 2 is required.  

{\bf Lemma 3}. Notations are the same as before. Denote
$g(x)=f(x)H_k(\frac{t-x}{\sqrt{2}\sigma})$, if $g(x)$ satisfies 
\begin{equation}
g'(x)\geq g(x)\frac{(x-t)}{{\sigma}^2}, 
\end{equation} 
whenever $x\leq t-M_2\Delta$,
then $\{F_k(l)\}_{l\in \mathcal{N}}$ is an increasing sequence.

The proof is also direct.  Therefore, by the same treatment as 
that in the previous subsection, we obtain
\begin{equation} \label{eq41}
\| E_3(t)\|_{L^2(R)}\leq \frac{\| f(t)\|_{L^2(R)}\sum_{i+j+k=s}
\frac{s!{\pi}^{i-1}H_k(\frac{-M_2\Delta}{\sqrt{2}\sigma})}
{i!k!{\Delta}^{i-1}(\sqrt{2}\sigma)^k(M_2\Delta)^{j+1}}}
{\exp\left(\frac{(M_2\Delta)^2}{2{\sigma}^2}\right)}. 
\end{equation}

\subsection{The end of the proof}

By combining Eqs. (\ref{eq10}), (\ref{eq25}), 
(\ref{eq39}) and (\ref{eq41}),  one obtains Eq. (\ref{eq3}). 

\section{Discussion}

\emph{Remark 1}.
For $s=0$, one has
\begin{eqnarray} 
&&\|f(t)-\sum_{n=\lceil\frac{t}{\Delta}\rceil-M_2}^{\lceil\frac{t}
{\Delta}\rceil+M_1}f(n\Delta)
\frac{\sin\left(\frac{\pi}{\Delta}(t-n\Delta)\right)}{\frac{\pi}
{\Delta}(t-n\Delta)}
\exp\left(-\frac{(t-n\Delta)^2}{2{\sigma}^2}\right)\|_{L^2(R)}
\nonumber\\ 
&&\leq
\sqrt{3}\| f(t)\|_{L^2(R)}\left\{\frac{1}{2\pi\sigma(\frac{\pi}
{\Delta}-B)
\exp\left(\frac{{\sigma}^2(\frac{\pi}{\Delta}-B)^2}{2}\right)}\right.
\nonumber\\
&&\left.+ \frac{1}{(M_1-1)\pi \exp\left(\frac{(M_1\Delta)^2}{2{
\sigma}^2}\right)}
+\frac{1}{M_2\pi \exp\left(\frac{(M_2\Delta)^2}{2{\sigma}^2}\right)}
\right\}.
\end{eqnarray} 
This is a rigorous error statement for 
the formulae widely used in the aforementioned  numerical 
computations.  
Roughly speaking, if $\exp(-\frac{x^2}{2})=10^{-\eta}$, 
then $\eta=\frac{x^2}{2ln10}$, so the error is 
\begin{eqnarray}
\varepsilon(r,B,\Delta,M)=\sqrt{3}\| f(t)\|_{L^2(R)}\left(\frac{1}
{2\pi r(\pi-B\Delta)
10^{\frac{r^2(\pi-B\Delta)^2}{2ln10}}}\right.\nonumber\\ 
\left.+ 
\frac{1}{(M_1-1)\Delta 10^{\frac{(M_1-1)^2}{2r^2ln10}}}
+\frac{1}{M_2\Delta10^{\frac{(M_2)^2}{2r^2ln10}}}\right), 
\end{eqnarray} 
where
$r=\frac{\sigma}{\Delta}$.  One may choose $r,B,\Delta$ and $ M$ 
appropriately to attain desired accuracy. 
Assume all non-exponential quantities are 
combined to give unit, and $M=M_1-1=M_2$, one has 
\begin{equation}\label{rule1} 
r(\pi-B\Delta)>\sqrt{\eta2\ln10}
\end{equation} 
and
\begin{equation} \label{rule2} 
\frac{M}{r}>\sqrt{\eta2\ln10},
\end{equation} 
where $\eta$ is the desired order of accuracy. 
There are some general rules for attaining high accuracy. 
These are discussed from two different arguments.

1.) For a given function $f(x)$ with a known bandlimit $B$,
      other parameters, $\Delta, r$ and $M$,
      are to be chosen appropriately 
      to achieve a desired accuracy
      order $\eta$:

(i) From Eqs. (\ref{rule1}) and (\ref{rule2}) one has 
   $B\Delta\leq \pi-\frac{\sqrt{\eta2\ln10}}{r}$. 
   For fixed $r$, the higher the frequency 
   bandlimit $B$ is, the smaller $\Delta$ should be, which means
   the more grid points in the computational domain.
   When $\Delta$ varies from $0$ to $\frac{\pi}{B}$, 
   $r$ changes from $\frac{\sqrt{\eta2\ln10}}{\pi}$ to $+\infty$, 
   therefore for sufficiently small $\Delta$, $r$ is 
   near $\frac{\sqrt{\eta2\ln10}}{\pi}$.

(ii) No matter how many grid points are in the computational domain, 
    $r$ and $M$ cannot be too small. Equations (\ref{rule1}) 
     and (\ref{rule2}) indicate
    $r>\frac{\sqrt{\eta2\ln10}}{\pi}$ and $M>r\sqrt{\eta2\ln10}$.  
    If $M$ and $r$ are less than the minimal requirements, the
    accuracy deteriorates quickly. 
    On the other hand, if sufficiently large  $r$ and $M$ are 
    used, say, $M=30 $ and $ r=3.5$, high approximation accuracy
    can be achieved.

2.) In practical computations, such as in solving a partial differential
    equation, the function $f$ and its $B$ are unknown.
    In this case, $\Delta$ is selected a priori. Then  $r$ and $M$
    are to be chosen properly for achieving 
    a desired accuracy order $\eta$:

(i)  For a given grid spacing, $\Delta$, and 
      accuracy requirement $\eta$,  $r$ value 
     determines frequency bandlimit $B$ which can be reached.
     Then the set of functions $f$ which are almost
     bandlimited to $B$ can be accurately approximated 
     (where `almost bandlimited to $B$' means the function $f$ 
     is not necessarily bandlimited but its Fourier
      amplitude outside $|B|$ is much smaller than the given 
     error $10^{-\eta}$).
     The choice of $M$ should be consistent with $r$ for 
     a given accuracy requirement. In general, 
     small $r$ and $M$
     values lead to an accurate approximation for 
     low frequency component of a function of interest.
     But the prediction of a high frequency component
     will not be accurate in such a case.

(ii) For a given grid spacing $\Delta$ and 
      $r$ value, the larger $M$ is, the higher
       bandlimit $B$ can be reached.

(iii)  To improve computational efficiency with a given $\Delta$, 
      $B$ shall be very close to 
     ${\pi\over\Delta}$. However, to maintain certain 
     approximation accuracy, $r$ 
     has to be sufficiently large, which implies that 
     $M$ has to be very large too. This in turn results in low 
     efficiency (It takes $M\rightarrow\infty$ 
     to maintain the accuracy if one samples at 
     the Nyquist rate).

(iv) If $\Delta$, $M$ and $\eta$ are chosen, then $r$ is also 
      fixed. For example, to achieve the machine 
      precision $10^{-\eta}\sim 10^{-15}$,   
      Eq. (\ref{rule1}) estimates $r>2.8$. If this is achieved
      by using $M=33$, then Eq. (\ref{rule2}) estimates 
      $ r< 4$. In fact, $M \sim 30$ and $2.8< r <4.0$ are 
      the parameter regions  found from an earlier numerical 
      test\cite{weicpl} and were used in many 
      applications\cite{weiphysica,wanwei,weigu,weijpb,GLW}.

\emph{Remark 2}. 
A comparison between the truncation 
errors of Shannon's sampling formulae and the regularized 
Shannon's sampling formulae is in order.
Reference \cite{Marks} estimates that the expression
\begin{equation}
(T_Nf)(t)=f(t)-\sum_{n=-N}^{n=+N}f(n\Delta)\frac{\sin(\frac{\pi}
{\Delta}(t-n\Delta))}
{\frac{\pi}{\Delta}(t-n\Delta)} 
\end{equation} 
has error of 
\begin{equation}
|T_N(t)|\leq \frac{\sqrt{2}}{\pi}\sqrt{E}\left|\sin(\frac{\pi
t}{\Delta})\right|\sqrt{\frac{N\Delta} {({N\Delta}^2-t^2)}}, 
\end{equation} 
where
$t<N\Delta$, and $E$ is the total `energy' of the function given by 
\begin{equation}
E=\int_{-\frac{\pi}{\Delta}}^{\frac{\pi}{\Delta}}|f\hat(w)|^2dw. 
\end{equation} 
This is not {\it directly} comparable with our error estimate because 
our sampling is centered around a point of interest, $x$. 
Let consider a truncation error of the form
\begin{equation}
(E_Mf)(t)=f(t)-\sum_{n=\lceil\frac{t}{\Delta}\rceil-M}^{n=\lceil\frac{t}
{\Delta}\rceil+M}
f(n\Delta)\frac{\sin(\frac{\pi}{\Delta}(t-n\Delta))} {\frac{\pi}
{\Delta}(t-n\Delta)}.
\end{equation} 
In Appendix A, it is shown that in a finite
computational domain, the $L^2 $ norm of $(E_Mf)(t)$ has 
the order of $\|f(t)\|_{L^2(R)} \sqrt{\frac{1}{M\Delta}}$, 
which is much larger than the
truncation error of the regularized Shannon's formulae.  
On the other hand, to achieve the same accuracy, the 
regularized formulae requires much fewer computational
grids\cite{weijcp,weicpl}.

\emph{Remark 3}. 
Discussions for the higher order derivatives 
can be presented in a similar manner as those of 
Remarks 1 and 2. In fact, previous 
work of solving partial differential 
equations\cite{weijcp,weicpl,weiphysica,wanwei,weigu,weijpb,GLW}.
involved such derivatives,
and results are consistent with 
the present theorem. Detailed comparison is 
omitted. 

\emph{Remark 4}. 
In many practical applications, such as in solving 
partial differential equations, error estimations 
and discussions
in other spaces are often required. Moreover, in 
real computations, the computational domain is 
always limited to a finite interval, such as $[a,b]$.
Therefore, the norm 
$\| f\|_{L^2(R)}$ in Eqs. (\ref{eq39}) and  (\ref{eq41}) are required 
to be changed into $\| f\|_{L^2(a,b)}$, which can be evaluated by 
integrations along $[a+M_1\Delta,b+M_1\Delta]$ and
$[a-M_2\Delta,b-M_2\Delta]$ respectively.  
Therefore various
$L^p (1\leq p\leq 2)$ error estimates of $E(t)$ 
can be derived accordingly.
If we know the size of $L^p$ norm 
$(1\leq p\leq 2)$ of the function of interest, then we can deduce 
from the theorem the critical values, $r$ and  $M$, to 
achieve desired accuracy.


\appendix

\section{Truncation error of Shannon's sampling formulae}

{\bf Lemma 4}. In the computational domain $[a,b]$, 
the error expression
\begin{equation}
(E_Mf)(t)=f(t)-\sum_{n=\lceil\frac{t}{\Delta}
\rceil-M}^{n=\lceil\frac{t}
{\Delta}\rceil+M}
f(n\Delta)\frac{\sin(\frac{\pi}{\Delta}(t-n\Delta))} 
{\frac{\pi}{\Delta}(t-n\Delta)}
\end{equation} 
satisfies the estimate 
\begin{equation} 
\|(E_Mf)(t)\|_{L^2(a,b)}\leq \frac{2\|
f(t)\|_{L^2(R)}}{\sqrt{(M-2)\Delta}}. 
\end{equation} 

{\bf Proof}. Let denote
\begin{equation} \label{eq52}
f_M(t)=\frac{1}{\pi}\sin\left(\frac{\pi t}
{\Delta}\right)
\sum_{n=\lceil\frac{t}{\Delta}\rceil-M}^{n=\lceil\frac{t}
{\Delta}\rceil+M}
\frac{f(n\Delta)(-1)^n}{(\frac{t}{\Delta}-n)}.
\end{equation} 
By using Schwartz's inequality, one has 
\begin{eqnarray}
&&\left(\sum_{n=\lceil\frac{t}{\Delta}\rceil+M}^{n=+\infty}
\frac{f(n\Delta)(-1)^n}
{(\frac{t}{\Delta}-n)}\right)^2\nonumber\\
&&\leq \|
f(t)\|_{L^2(R)}^2\sum_{n=\lceil\frac{t}{\Delta}
\rceil+M}^{n=+\infty}\frac{1}
{(t-n\Delta)^2}\nonumber\\
&&\leq \frac{1}{{\Delta}^2}\| f(t)\|_{L^2(R)}^2\sum_{l\geq
M}\frac{1}{(l-1)^2}\nonumber\\ 
&&\leq\frac{1}{{\Delta}^2}\|
f(t)\|_{L^2(R)}^2\int_{M-2}^{+\infty}\frac{dx}{x^2}\nonumber\\ 
&&=\frac{\|
f(t)\|_{L^2(R)}^2}{(M-2){\Delta}^2}. \label{eq53}
\end{eqnarray} 
Similarly one obtains
\begin{equation}\label{eq54}
\left(\sum_{n=-\infty}^{n=\lceil\frac{t}{\Delta}\rceil-M}
\frac{f(n\Delta)(-1)^n}
{(\frac{t}{\Delta}-n)}\right)^2
\leq\frac{\| f(t)\|_{L^2(R)}^2}{(M-1){\Delta}^2}. 
\end{equation} 
By combining Eqs. (\ref{eq52}), (\ref{eq53}) and
(\ref{eq54}), one finishes the proof. 



\begin{thebibliography}{99} 
\bibitem{weijcp}G. W. Wei, ``Discrete singular convolution
for the solution of the Fokker-Planck equation'', \emph{J. Chem.  
   Phys}.,  vol. 110, 8930-8942 (1999). 

\bibitem{Shannon}C. E. Shannon, ``A mathematical theory of communication'',
\emph{Bell System Tech. J.,} 27, 379-423 (1948). 

\bibitem{Chui}Charles K. Chui, \emph{An Introduction to 
           Wavelets}, Academic Press, 1992.


\bibitem{Marks}Robert J. Marks II, \emph{Introduction to Shannon 
        Sampling and Interpolation Theory}, Spring-Verlag, 1991. 


\bibitem{weiprl} G. W. Wei, D. S. Zhang, D. J. Kouri and
    D. K. Hoffman, ``Lagrange Distributed Approximating Functionals'', 
  \emph{Phys. Rev. Lett.}, vol. 79, 775-779 (1997). 


\bibitem{weicpl}G. W. Wei, ``Quasi wavelets and quasi interpolating 
   wavelets'', \emph{Chem. Phys.  Lett}., vol. 296, 215-222 (1998). 


\bibitem{weiphysica}G. W. Wei, 
          ``Discrete singular convolution method for the Sine-Gordon
          equation'',         
        { \it  Physica} D, vol. 137, 247-259 (2000).         
  

\bibitem{wanwei}G. W. Wei 
        ``A unified computational method for solving mechanical problems'',
            { \it Commun. Numer. Methods Engng.}, submitted;         
       D. C. Wan and G. W. Wei, 
       ``Numerical solution of unsteady incompressible flows by the 
       discrete singular convolution'',
       {\it Int. J. Numer. Methods Fluid}, submitted.         


\bibitem{weigu}G. W. Wei and Y. Gu, ``A novel approach for Burgers' equation 
        with high Reynolds number'', unpublished.
    


\bibitem{weijpb}G. W. Wei, 
          ``Solving quantum eigenvalue problems by 
           discrete singular convolution'',
            { \it J. Phys. B}, vol. 33, 343-352 (2000).         


\bibitem{GLW}S. Guan, C.-H. Lai and G. W. Wei, 
        ``Boundary Controlled Nanoscale Morphology
       in a Circular Domain'',
        { \it Physica D}, submitted.
 
\bibitem{Zwillinger}D. Zwillinger,  \emph{Handbook of Integration}, 
Jones and Barlett, 1992. 

\bibitem{Mitrinovic}D. S. Mitrinovic, \emph{Analytic Inequalities}, 
Spring-Verlag, 1970.



\end{thebibliography}
\end{document}